\documentclass[12pt,thmsa]{article}
\usepackage{sw20lart}

\input{tcilatex}
\begin{document}

\author{Alan Horwitz}
\title{A Version of Simpson's Rule for Multiple Integrals}
\maketitle

\section{Introduction}

\footnotetext{
Key Words: Cubature Rule, Simpson's Rule, polygonal region, Grobner Basis,
exact}Let $M(f)$ denote the Midpoint Rule and $T(f)$ the Trapezoidal Rule
for estimating $\int_{a}^{b}f(x)dx$. Then $M(f)$ and $T(f)$ are each exact
for $f(x)=1$ and $x$. A more accurate rule can be obtained by taking $%
\lambda M(f)+(1-\lambda )T(f)$, where $\lambda =\frac{2}{3}$. Of course this
rule is known as Simpson's Rule, and is exact for all polynomials of degree $%
\leq 3$. The purpose of this paper is to extend this idea to multiple
integrals over certain polygonal regions $D_{n}$ in $R^{n}$. First, consider
the case $n=2$. So suppose that $D$ is a polygonal region in $R^{2}$. The
midpoint rule in one dimension is $M(f)=(b-a)f\left( \frac{a+b}{2}\right) .$
A natural extension of this rule is 
\[
M(f)=A(D)f(P) 
\]
where $A(D)=$ area of $D$, and $P=$ centroid of $D$.

The trapezoidal rule in one dimension is $T(f)=(b-a)\left( \frac{f(a)+f(b)}{2%
}\right) .$ A natural extension of that rule is 
\[
T(f)=A(D)\frac{1}{m}\sum_{k=1}^{m}f(P_{k}) 
\]
where the $P_{k}$ are the vertices of $D_{n}$.

In general, let $D_{n}$ be some polygonal region in $R^{n},$ let $%
P_{0},...,P_{m}$ denote the vertices of $D_{n}$, and let $P_{m+1}=$ center
of mass of $D_{n}$. Define the linear functionals

$M(f)=$ Vol$(D_{n})$ $f(P_{m+1}),$ $T(f)=$ Vol$(D_{n})(\frac{1}{m+1}%
\sum_{j=0}^{m}f(P_{j})),$ and for fixed $\lambda ,\;0\leq \lambda \leq 1$,

\[
L_{\lambda }=\lambda M(f)+(1-\lambda )T(f) 
\]
The idea is to choose $\lambda $ so that $L_{\lambda }$ is a good cubature
rule(CR). Our objective here is not to provide \textbf{optimal} cubature
rules, but to generalize the ideas from Simpson's Rule in one variable to
several variables. In some cases we have reproduced known CRs using a
different approach\footnote{%
See \cite{2} and \cite{5} for a thorough treatment of CRs.}. In other cases
our CRs appear to be new. In addition, our approach suggest a general method
for deriving CRs for polygonal domains.

For the $n$ cube $[0,1]\times \cdots \times [0,1]$, $\lambda =\frac{2}{3}$,
and the CR $L(f)$ is exact for all polynomials of degree $\leq 3$(see CR3).
This, of course, matches the degree of exactness of Simpson's Rule in one
variable. For the $n$ simplex, $\lambda =\frac{n+1}{n+2}$, and the CR $L(f)$
is exact for all polynomials of degree $\leq 2$(see CR1).

For regular polygons in the plane with $m$ sides, our approach fails for $%
m>4 $. For example, if $D=$ a regular hexagon, then \textbf{no} choice of $%
\lambda $ makes $L_{\lambda }$ exact for all polynomials of degree $\leq 2$.
Even for $m=4$, if the polygon is \textbf{not regular}, then our approach
fails as well. For example, let $D$ be the trapezoid with vertices $%
\{(0,0),(1,0),(0,1),(1,2)\}$. Again, \textbf{no} choice of $\lambda $ makes $%
L_{\lambda }$ exact for all polynomials of degree $\leq 2$. One can do
better, however, by using points $Q_{j}\in \partial (D_{n}),$ \textbf{other
than the vertices} of $D_{n}$, to generate $T(f)$. In some cases this leads
to better formulas. For example, for the trapezoid $D$, this leads to a
formula which is exact for polynomials of degree $\leq 2$(see CR 5)$.$
Choosing points different from the vertices leads to a system of polynomial
equations. We use Grobner Basis methods(see \cite{3}), along with Maple, to
solve such systems when possible. For the $n$ Simplex $T_{n}$, we try using
the center of mass of the faces of $T_{n}$ to generate $T(f)$(see CR2). A
similar idea works for the $n$ cube.

Instead of using the weighted combination, $\lambda M(f)+(1-\lambda )T(f)$,
another way to derive Simpson's Rule in one variable is to integrate the
quadratic interpolant to $f$ at $a,\frac{a+b}{2},$ and $b.$ For the $n$
Simplex $T_{n}$, our generalization $L_{\lambda }$ can be obtained by
integrating $p(\widehat{x})$ over $T_{n}$, where $p(\widehat{x})$ is the
unique interpolant to $f(\widehat{x})$ at the vertices and center of mass of 
$T_{n}$. Indeed, $M(f)=\int_{T_{n}}T(\widehat{x})dV$, where $T(\widehat{x})=$
tangent hyperplane to $f(\widehat{x})$ at the center of mass of $T_{n}$, and 
$T(f)=\int_{T_{n}}L(\widehat{x})dV$, where $L(\widehat{x})=$ bilinear
interpolant to $f(\widehat{x})$ at the vertices of $T_{n}$. This, of course,
is the generalization from one variable, where $(b-a)f\left( \frac{a+b}{2}%
\right) =\int_{a}^{b}T(x)dx$, $T=$ tangent line to $f(x)$ at $\frac{a+b}{2}$%
, and the trapezoidal rule can be obtained by integrating the linear
interpolant to $f$ at the endpoints. For regions $D_{n}$ in general,
however, $L_{\lambda },T,$ and $M$ do not always arise in this fashion.
Indeed, the number of basis functions does not always match the number of
nodes, and the interpolant may not be unique. This is what occurs for the $n$
cube.

We also give an analgous formula, using four knots, for the unit disc(see CR
6).

Finally we note that for most of our CRs, all of the weights are positive,
all of the knots lie inside the region $D_{n}$, and all but one of the knots
lies on the boundary of the region.

\section{$n$ Simplex}

Let $P_{0},...,P_{n}$ denote the $n+1$ vertices of the $n$ simplex $%
T_{n}\subset R^{n}$. Hence $P_{j}=(0,0,...,1,0,0)$ for $\dot{k}\geq 1$, and $%
P_{0}=(0,0,...,0),\widehat{x}=(x_{1},...,x_{n}),dV=$ standard Lebesgue
measure on $T_{n}$. We list the following useful facts:

\[
\text{Vol}(T_{n})=\frac{1}{n!}\text{, }P_{n+1}=\text{Center of Mass of }%
T_{n}=(1/(n+1),...,1/(n+1)) 
\]
\[
\int_{T_{n}}x_{k}dV=\frac{1}{(n+1)!},k=1,...,n 
\]
\[
\int_{T_{n}}x_{j}x_{k}dV=\frac{1}{(n+2)!},j\neq k 
\]
\[
\int_{T_{n}}x_{k}^{2}dV=\frac{2}{(n+2)!},k=1,...,n 
\]
Define the following linear functionals:

\[
M(f)=\text{Vol}(T_{n})\text{ }f(P_{n+1})=\frac{1}{n!}f(1/(n+1),...,1/(n+1)) 
\]
\[
T(f)=\text{Vol}(T_{n})(\frac{1}{n+1}\sum_{j=0}^{n}f(P_{j}))=\frac{1}{(n+1)!}%
\sum_{j=0}^{n}f(P_{j}) 
\]
\[
I(f)=\int_{T_{n}}f(\widehat{x})dV 
\]
Finally, for fixed $\lambda ,\;0\leq \lambda \leq 1$, define

\[
L_{\lambda }=\lambda M(f)+(1-\lambda )T(f) 
\]
First, if $f(\widehat{x})=x_{k}$, then $I(f)=M(f)=T(f)=\frac{1}{(n+1)!}$ for
any $\lambda $. Now let $f(\widehat{x})=x_{j}x_{k},\;j\neq k$. Then $I(f)=%
\frac{1}{(n+2)!},M(f)=\frac{1}{(n+1)(n+1)!},T(f)=0$. Hence $L_{\lambda
}(f)=I(f)\Rightarrow $ $\lambda \left( \frac{1}{(n+1)(n+1)!}\right) =\frac{1%
}{(n+2)!}\Rightarrow $ $\lambda =\frac{n+1}{n+2}$. Let $L=L_{\lambda },$with
this $\lambda .$ Now, if $f(\widehat{x})=x_{k}^{2}$, then $L(f)=\frac{n+1}{%
n+2}\frac{1}{n!}\frac{1}{(n+1)^{2}}+\frac{1}{n+2}\frac{1}{(n+1)!}=\frac{2}{%
(n+2)!}=I(f)$. We summarize

\[
\text{\textbf{CR1:} }L(f)=\frac{n+1}{(n+2)n!}f(P_{n+1})+\frac{1}{(n+2)!}%
\sum_{j=0}^{n}f(P_{j})\text{ is exact for all polynomials of degree}\leq 2. 
\]
$L$ is \textit{not} exact for all polynomials of degree $\leq 3$. For
example, if $f(x)=x_{j}x_{k}x_{l},j\neq k,k\neq l,j\neq l,$ then $I(f)=\frac{%
1}{(n+3)!}$, while $L(f)=\frac{1}{(n+1)(n+2)!}\neq I(f)$.

\subsubsection{Connection with Interpolation}

Another way to derive Simpson's Rule in one variable is by using quadratic
interpolation at $a,\frac{a+b}{2},b$. The cubature rule $L(f)$ above can
also be obtained by integrating a certain second degree interpolant to $f$
over $T_{n}$. Let 
\[
p(x_{1},...,x_{n})=A_{n+1}x_{1}x_{2}+\sum_{k=1}^{n}A_{k}x_{k}+A_{0} 
\]
We wish to choose the $A_{j}$ so that 
\[
p(P_{k})=f(P_{k}),\;k=0,...,n+1 
\]
\label{eqS1}

It follows easily from (\ref{eqS1}) that $%
A_{0}=f(P_{0}),A_{k}=f(P_{k}),k=1,...,n,$ and

$A_{n+1}=f(P_{n+1})-(n+1)\sum_{k=0}^{n}f(P_{k})$. Hence $\int_{T_{n}}p(%
\widehat{x})dV=A_{n+2}\frac{1}{(n+2)!}+\frac{1}{(n+1)!}\sum_{k=1}^{n}A_{k}+%
\frac{1}{n!}A_{0}=\frac{1}{(n+2)!}\left(
(n+1)^{2}f(P_{n+1})-(n+1)\sum_{k=0}^{n}f(P_{k})\right) $

$+\frac{1}{(n+1)!}\left( \sum_{k=1}^{n}f(P_{k})-nf(P_{0})\right) +\frac{1}{n!%
}f(P_{0})=\frac{(n+1)^{2}}{(n+2)!}f(P_{n+1})+\frac{1}{(n+2)!}%
\sum_{k=0}^{n}f(P_{k})=L(f)$.

\begin{remark}
Note that one could use any $x_{j}x_{k},j\neq k$ in constructing the
interpolant $p,$ and still obtain $\int_{T_{n}}p(\widehat{x})dV=L(f).$
\end{remark}

\begin{remark}
This, of course, generalizes the fact that Simpson's Rule in one variable
can be obtained by integrating the quadratic interpolant to $f$ at $a,b,%
\frac{a+b}{2}$. In more than one variable, however, our quadratic
interpolant only uses one second degree basis function. Using the basis
functions $\{1,x_{1},...,x_{n},x_{1}x_{2}\}$, one obtains exactness for
those functions as well as for $\frac{n^{2}+n-2}{2}$ additional quadratic
terms.
\end{remark}

\begin{remark}
One can also obtain the rules $M(f)$ and $T(f)$ using interpolation.
\end{remark}

$M(f)=\int_{T_{n}}T(\widehat{x})dV$, where $T$ is the tangent hyperplane to $%
f$ at $P_{n+1}$.

$T(f)=\int_{T_{n}}q(\widehat{x})dV$, where $q(x_{1},...,x_{n})=%
\sum_{k=1}^{n}B_{k}x_{k}+B_{0}$, the $B_{k}$ chosen so that

$q(P_{k})=f(P_{k}),\;k=0,...,n$.

Again, this generalizes the fact that the midpoint rule in one variable can
be obtained by integrating the tangent line at the midpoint, while the
trapezoidal rule can be obtained by integrating the linear interpolant to $f$
at the endpoints.

\subsection{Other Boundary Points}

It is interesting to examine what happens if we use other points on $%
\partial (S_{n})$ to generate $T(f)$. We first examine $n=2$ in some detail,
and then generalize.

\subsubsection{$n=2$}

Let $T_{2}=$ the triangle in $R^{2}$ with vertices $\{(0,0),(1,0),(0,1)\}$.
Use the points $(a,0),(0,b),(c,1-c),0\leq a\leq 1,0\leq b\leq 1,0\leq c\leq
1 $ from each side of $T_{2}$ to define

\[
T(f)=\text{Area}(T_{2})\frac{1}{3}(f(a,0)+f(0,b)+f(c,1-c))=\frac{1}{6}%
(f(a,0)+f(0,b)+f(c,1-c)) 
\]
As earlier(with $(1/3,1/3)=$ Center of mass of $T$), 
\[
M(f)=\text{Area}(T_{2})f(1/3,1/3)=\frac{1}{2}f(1/3,1/3) 
\]
and, for fixed $\lambda ,\;0\leq \lambda \leq 1$, 
\[
L_{\lambda }=\lambda M(f)+(1-\lambda )T(f) 
\]

\begin{itemize}
\item  $L_{\lambda }(x)=\allowbreak \frac{1}{6}\lambda +\frac{1}{6}\left(
1-\lambda \right) \left( a+c\right) $

\item  $L_{\lambda }(y)=\frac{1}{6}\lambda +\frac{1}{6}\left( 1-\lambda
\right) \left( b+1-c\right) $

\item  $\allowbreak L_{\lambda }(x^{2})=\allowbreak \frac{1}{18}\lambda +%
\frac{1}{6}\left( 1-\lambda \right) \left( a^{2}+c^{2}\right) $

\item  $L_{\lambda }(y^{2})=\allowbreak \frac{1}{18}\lambda +\frac{1}{6}%
\left( 1-\lambda \right) \left( b^{2}+\left( 1-c\right) ^{2}\right) $

\item  $L_{\lambda }(xy)=\allowbreak \frac{1}{18}\lambda +\frac{1}{6}\left(
1-\lambda \right) c\left( 1-c\right) $
\end{itemize}

We want to choose $a,b,c,$ and $\lambda $ so that $L_{\lambda }$ is exact
for the functions $x,y,x^{2},$ $y^{2}$,and $xy$. Setting $L_{\lambda
}(f)=I(f)$ yields the system of equations

$\allowbreak \frac{1}{6}\lambda +\frac{1}{6}\left( 1-\lambda \right) \left(
a+c\right) =\frac{1}{6}$, $\frac{1}{6}\lambda +\frac{1}{6}\left( 1-\lambda
\right) \left( b+1-c\right) =\frac{1}{6}$, $\allowbreak $

$\frac{1}{18}\lambda +\frac{1}{6}\left( 1-\lambda \right) \left(
a^{2}+c^{2}\right) =\frac{1}{12}$, $\frac{1}{18}\lambda +\frac{1}{6}\left(
1-\lambda \right) \left( b^{2}+\left( 1-c\right) ^{2}\right) =\frac{1}{12}$,

$\frac{1}{18}\lambda +\frac{1}{6}\left( 1-\lambda \right) c\left( 1-c\right)
=\frac{1}{24}$

We found the following Grobner basis using Maple:

\[
\{-12c^{2}+12\lambda c^{2}+4\lambda -3+12c-12\lambda c,-1+a+c,b-c\} 
\]
Setting each polynomial from the Grobner basis to $0$ yields precisely the
same solutions as the original system.

Hence $b=c$, $a=1-c$, and $\lambda =\frac{12c^{2}+3-12c}{-12c^{2}+4-12c}%
\allowbreak $, $1-\lambda =\allowbreak \frac{1}{4}\frac{24c^{2}-1}{%
3c^{2}-1+3c}$

This gives $L_{\lambda }(f)=\lambda \frac{1}{2}f(1/3,1/3)+(1-\lambda )\frac{1%
}{6}(f(1-c,0)+f(0,c)+f(c,1-c))$

Now we choose $c$ so that $L_{\lambda }$ is exact for quadratics. Since $%
L_{\lambda }(x^{2})=\allowbreak \frac{1}{12}\frac{9c^{2}-1+3c-24c^{3}+24c^{4}%
}{3c^{2}-1+3c}$ and $\int_{T_{2}}x^{2}dA=\allowbreak \frac{1}{12}$, $%
L_{\lambda }(x^{2})=\frac{1}{12}\Rightarrow c=0$ or $c=\frac{1}{2}$. $c=0$
uses the vertices of $T_{2}$ to generate $T(f)$. Hence we use $c=\frac{1}{2}%
\Rightarrow \lambda =\allowbreak 0$. With this choice of $c$ and $\lambda $,
it follows that $L_{\lambda }(y^{2})=\int_{T_{2}}y^{2}dA$ and $L_{\lambda
}(xy)=\int_{T_{2}}xydA$.

We summarize

\[
L(f)=\frac{1}{6}(f(\frac{1}{2},0)+f(0,\frac{1}{2})+f(\frac{1}{2},\frac{1}{2}%
))\text{ is exact for all polynomials of degree}\leq 2. 
\]

\begin{remark}
The CR above is given in the well known book of Stroud(see \cite{5}), as one
of a group of formulas for the triangle $T_{2}$.
\end{remark}

Since $L(x^{3})=\frac{1}{24}\neq \int_{T_{2}}x^{3}dA=\allowbreak \frac{1}{20}
$, $L$ is not exact for all cubics.

\paragraph{Connection with Interpolation}

It is not hard to show that there is a unique interpolant, $%
p(x,y)=A_{1}xy+A_{2}x+A_{3}y+A_{4}$, to any given $f(x,y)$ at $%
(1/2,0),(0,1/2),(1/2,1/2),(1/3,1/3)$. It also follows easily that $%
\int_{T_{2}}p(x,y)dA=\frac{1}{6}(f(\frac{1}{2},0)+f(0,\frac{1}{2})+f(\frac{1%
}{2},\frac{1}{2}))$, which is $L(f)$ given above.

\subsubsection{General $n$}

For $T_{2}$ above we used the midpoints of the edges for $T(f)$, so it is
natural to try using the \textbf{center of mass of the faces} of $T_{n}$,
given by $Q_{k}=(\frac{1}{n},\frac{1}{n},...,0,...,\frac{1}{n})$(all
coordinates $\frac{1}{n}$ except the $k$th coordinate $=0$), $k=1,...,n$, $%
Q_{n+1}=(\frac{1}{n},\frac{1}{n},...,\frac{1}{n}).$ This yields

\[
T(f)=\frac{1}{(n+1)!}\sum_{k=1}^{n+1}f(Q_{k}) 
\]
As earlier,

\[
M(f)=\frac{1}{n!}f(1/(n+1),...,1/(n+1)) 
\]
\[
L_{\lambda }=\lambda M(f)+(1-\lambda )T(f) 
\]
A simple computation shows that $L_{\lambda }$ is exact for $x_{k}$ for any $%
\lambda ,$and is exact for $x_{j}x_{k},j\neq k$, if $(\frac{1}{n+1})\lambda
+(\frac{n-1}{n^{2}})(1-\lambda )=\frac{1}{n+2}.$ This gives $\lambda
=\allowbreak -\left( n-2\right) \frac{n+1}{n+2}$. With this $\lambda $, we
have

\[
\text{\textbf{CR2}: }L(f)=-\left( n-2\right) \frac{n+1}{n+2}\frac{1}{n!}%
f(1/(n+1),...,1/(n+1))+\frac{n^{2}}{(n+2)!}\sum_{k=1}^{n+1}f(Q_{k})\text{ } 
\]
is exact for all polynomials of degree $\leq 2.$ Note that one of the
weights is $<0$ if $n\geq 3$. $\allowbreak $

\begin{remark}
It would be nice to choose points on the faces which yield $\lambda \geq 0$,
and in particular $\lambda =0$. We tried this for $n=3$, using points on
each face of the simplex $T_{3}$. Letting $%
Q_{1}=(a_{1},a_{2},0),Q_{2}=(a_{3},0,a_{4}),Q_{3}=(0,a_{5},a_{6}),Q_{4}=(a_{7},a_{8},1-a_{7}-a_{8}) 
$, with
\end{remark}

\[
a_{j}\geq 0\;\forall j,\;a_{j}+a_{j+1}\leq 1,j=1,3,5,7 
\]
\[
T(f)=\text{area}(T_{3})\frac{1}{4}\sum_{j=1}^{4}f(Q_{j})=\frac{1}{24}%
\sum_{j=1}^{4}f(Q_{j}) 
\]
\[
M(f)=\frac{1}{6}f(1/4,1/4,1/4) 
\]
\[
L_{\lambda }=\lambda M(f)+(1-\lambda )T(f) 
\]

The equations giving exactness for $x,y,z,x^{2},y^{2},z^{2},xy,xz,yz$, after
some simplification, look like $%
a_{1}+a_{3}+a_{7}=1,a_{2}+a_{5}+a_{8}=1,a_{4}+a_{6}-a_{7}-a_{8}=0,$

$\frac{1}{96}$ $\lambda +(1-$ $\lambda )\frac{1}{24}(a_{1}a_{2}+a_{7}a_{8})=%
\frac{1}{120},$

$\frac{1}{96}$ $\lambda +(1-\lambda )\frac{1}{24}%
(a_{3}a_{4}+a_{7}(1-a_{7}-a_{8}))=\frac{1}{120},$

$\frac{1}{96}$ $\lambda +(1-$ $\lambda )\frac{1}{24}%
(a_{5}a_{6}+a_{8}(1-a_{7}-a_{8}))=\frac{1}{120},$

$\frac{1}{96}$ $\lambda +(1-$ $\lambda )\frac{1}{24}%
(a_{1}^{2}+a_{3}^{2}+a_{7}^{2})=\frac{1}{60},$

$\frac{1}{96}$ $\lambda +(1-$ $\lambda )\frac{1}{24}%
(a_{2}^{2}+a_{5}^{2}+a_{8}^{2})=\frac{1}{60},$

$\frac{1}{96}$ $\lambda +(1-$ $\lambda )\frac{1}{24}%
(a_{4}^{2}+a_{6}^{2}+(1-a_{7}-a_{8})^{2})=\frac{1}{60}$

One can solve the first two equations for $a_{7}$ and $a_{8}$, and then
substitute into the rest of the equations, yielding seven polynomial
equations in seven unknowns. Of course two solutions to these equations are

\begin{itemize}
\item  $\;a_{1}=a_{2}=a_{4}=a_{6}=a_{7}=a_{8}=0,\;a_{3}=a_{4}=1,\;\lambda =%
\frac{4}{5}$, which uses the vertices of $T_{3}$ for $Q_{j}$ , and

\item  $\;a_{j}=\frac{1}{3}$ for all $j$, $,\;\lambda =-\frac{4}{5}$, which
uses the center of mass of each face of $T_{3}$ for $Q_{j}$.
\end{itemize}

We applied Grobner Basis methods to the above equations giving exactness for 
$x,y,z,x^{2},y^{2},z^{2},xy,xz,yz$, along with $\lambda =0$. We were then
able to show that there is no solution. That proves that one cannot make $%
\lambda =0$. We were not able to find any positive values for $\lambda $,
though it is possible that some exist.

\section{Unit $n$ Cube}

Let $C_{n}=[0,1]^{n}\subset R^{n}$, and let $P_{0},...,P_{m-1}$ denote the $%
m $ vertices of $C_{n},m=2^{n}$. Assume that $P_{0}=(0,...,0)$ and $%
P_{m-1}=(1,...,1)$. The other vertices have at least one coordinate which
equals $0$, and at least one coordinate which equals $1$. Let $\widehat{x}%
=(x_{1},...,x_{n}),dV=$ standard Lebesgue measure on $C_{n}$. We list the
following useful facts:

\[
\text{vol}(C_{n})=1,\;\text{Center of Mass}=P_{m}=(1/2,\ldots ,1/2) 
\]
For any non--negative $r_{1},...,r_{n}$,

\[
\int_{C_{n}}x_{k_{1}}^{r_{1}}x_{k_{2}}^{r_{2}}\cdots x_{k_{j}}^{r_{j}}dV=%
\frac{1}{r_{1}+1}\frac{1}{r_{2}+1}\cdots \frac{1}{r_{j}+1} 
\]
Define the following linear functionals:

\[
M(f)=\text{Vol}(C_{n})\text{ }f(P_{m})=f(1/2,...,1/2) 
\]
\[
T(f)=\text{Vol}(C_{n})(\frac{1}{m}\sum_{j=0}^{m-1}f(P_{j}))=\frac{1}{m}%
\sum_{j=0}^{m-1}f(P_{j}) 
\]
\[
I(f)=\int_{C_{n}}f(\widehat{x})dV 
\]
Finally, for fixed $\lambda ,\;0\leq \lambda \leq 1$, 
\[
L_{\lambda }=\lambda M(f)+(1-\lambda )T(f) 
\]
First, let $f(x)=x_{k}$. Then $M(f)=\frac{1}{2}$, and it is not hard to show
that $T(f)=\frac{1}{m}2^{n-1}=\frac{1}{2}$. Hence $M(f)=T(f)=I(f)\Rightarrow
L_{\lambda }(f)=I(f)$ for any $\lambda .$ Now let $f(x)=x_{j}x_{k},j\neq k$.
Then $M(f)=\frac{1}{4}$. since there are $2^{n-2}$ ways to get a $1$ in both
the $jth$ and $kth$ coordinates of a vertex of $D$, $T(f)=\frac{1}{m}2^{n-2}=%
\frac{1}{4}$. Hence, again, $M(f)=T(f)=I(f)$ for any $\lambda \Rightarrow
L_{\lambda }(f)=I(f).$ Also, if $f(x)=x_{k}^{2}$, then $M(f)=\frac{1}{4}$
and $T(f)=\frac{1}{2}$. Hence $L_{\lambda }(f)=I(f)=\frac{1}{3}\Rightarrow $ 
$\lambda =\frac{2}{3}$. Now let $L=L_{2/3}$. Finally we consider third
degree terms. First, if $f(x)=x_{j}x_{k}x_{l},j\neq k,k\neq l,j\neq l,$ then 
$T(f)=M(f)=I(f)=\frac{1}{8}\Rightarrow L(f)=I(f)$ for any $\lambda .$
Second, if $f(x)=x_{j}^{2}x_{k},j\neq k,$ then $T(f)=\frac{1}{4},M(f)=\frac{1%
}{8},$ and $I(f)=\frac{1}{6}$, and hence $\frac{2}{3}M(f)+\frac{1}{3}%
T(f)=I(f)$. Finally, if $f(x)=x_{k}^{3}$, then

$T(f)=\frac{1}{2},M(f)=\frac{1}{8},$ and $I(f)=\frac{1}{4}$, and again $%
\frac{2}{3}M(f)+\frac{1}{3}T(f)=I(f)$. We summarize

\[
\text{\textbf{CR3:} }L(f)=\frac{2}{3}f(1/2,...,1/2)+\frac{1}{3}\frac{1}{2^{n}%
}\sum_{j=0}^{2^{n}-1}f(P_{j})\text{ is exact for all polynomials of degree}%
\leq 3. 
\]

\begin{remark}
If $f(x)=x_{k}^{4}$, then $T(f)=\frac{1}{2},M(f)=\frac{1}{16},$ and $I(f)=%
\frac{1}{5}\Rightarrow \frac{2}{3}M(f)+\frac{1}{3}T(f)\neq I(f)$.
\end{remark}

Hence $L$ is \textit{not exact} in general for polynomials of degree $\leq 4$%
.

\subsection{Other Points on Boundary of $n$ Cube}

\subsubsection{$n=2$}

Let $D=$ unit square. The formulas above use the vertices of $D$ to generate 
$T(f)$. Here we use the points $(a,0),(0,b),(c,1),(1,d)$ on $\partial (D)$, $%
0\leq a\leq 1,0\leq b\leq 1,0\leq c\leq 1,0\leq d\leq 1,$ to define 
\[
T(f)=\frac{1}{4}(f(a,0)+f(0,b)+f(c,1)+f(1,d)) 
\]
As earlier,

\[
M(f)=f(1/2,1/2) 
\]
For$\;0\leq \lambda \leq 1$,

\[
L_{\lambda }=\lambda M(f)+(1-\lambda )T(f) 
\]
The idea is to choose $a,b,c,d$ and $\lambda $ so that $L_{\lambda }$ is
exact for all polynomials of degree $\leq 3$.

\begin{itemize}
\item  $\;L_{\lambda }(x)=\frac{1}{2}\lambda +\frac{1}{4}\left( 1-\lambda
\right) \left( a+c+1\right) $

\item  $\;L_{\lambda }(y)=\allowbreak \frac{1}{2}\lambda +\frac{1}{4}\left(
1-\lambda \right) \left( b+1+d\right) $

\item  $L_{\lambda }(x^{2})=\allowbreak \frac{1}{4}\lambda +\frac{1}{4}%
\left( 1-\lambda \right) \left( a^{2}+c^{2}+1\right) $

\item  $L_{\lambda }(y^{2})=\frac{1}{4}\lambda +\frac{1}{4}\left( 1-\lambda
\right) \left( b^{2}+1+d^{2}\right) $

\item  $L_{\lambda }(xy)=\allowbreak \frac{1}{4}\lambda +\frac{1}{4}\left(
1-\lambda \right) \left( c+d\right) $

\item  $L_{\lambda }(x^{3})=\allowbreak \frac{1}{8}\lambda +\frac{1}{4}%
\left( 1-\lambda \right) \left( a^{3}+c^{3}+1\right) $

\item  $L_{\lambda }(y^{3})=\allowbreak \frac{1}{8}\lambda +\frac{1}{4}%
\left( 1-\lambda \right) \left( b^{3}+1+d^{3}\right) $

\item  $L_{\lambda }(x^{2}y)=\allowbreak \frac{1}{8}\lambda +\frac{1}{4}%
\left( 1-\lambda \right) \left( c^{2}+d\right) $

\item  $L_{\lambda }(xy^{2})=\allowbreak \frac{1}{8}\lambda +\frac{1}{4}%
\left( 1-\lambda \right) \left( c+d^{2}\right) $
\end{itemize}

Letting $I(f)=\int_{0}^{1}\int_{0}^{1}f(x,y)dydx$, we have $%
I(x)=I(y)=\allowbreak \frac{1}{2}$, $I(x^{2})=I(y^{2})=\allowbreak \frac{1}{3%
}$, $I(xy)=I(x^{3})=\allowbreak I(y^{3})=\allowbreak \frac{1}{4}$, $%
I(x^{2}y)=\allowbreak I(xy^{2})=\frac{1}{6}$, $I(x^{2}y^{2})=\allowbreak 
\frac{1}{9}$, $I(x^{3}y)=\allowbreak I(xy^{3})=\frac{1}{8}$, $%
I(x^{4})=I(y^{4})=\allowbreak \frac{1}{5}$. Setting $L_{\lambda }(f)=I(f)$
for the nine monomials above yields a system of nine polynomial equations in
nine unknowns. Maple gives the Grobner Basis 
\[
\{3\lambda +6d-2-6\lambda d-6d^{2}+6\lambda d^{2},a-d,-1+b+d,-1+c+d\} 
\]
$a=d$, $b=1-d$, $c=1-d$, $\lambda =\frac{6d^{2}-6d+2}{6d^{2}-6d+3}$. We are
free to choose $d$ to force exactness for additional monomials. However, it
is not possible to get exactness for all fourth degree polynomials. Say we
want $L_{\lambda }(x^{3}y)=I(x^{3}y)=\frac{1}{8}\Rightarrow -\frac{1}{24}%
\frac{-9d^{2}+7d-3+2d^{3}}{2d^{2}-2d+1}=\frac{1}{8}\Rightarrow
d=0,\allowbreak d=1,\allowbreak d=\frac{1}{2}$. $d=0$ or $d=1$ uses the
vertices of $D$--the case $n=2$ of the general method discussed earlier for
the $n$ cube. For $\allowbreak d=\frac{1}{2}$ we have $\lambda =\allowbreak 
\frac{1}{3}$. with this choice, $L_{\lambda }$ is also exact for $%
f(x,y)=xy^{3}$. $L_{\lambda }(x^{4})=\allowbreak \frac{5}{24}\neq
I(x^{4})=\allowbreak \frac{1}{5}$, so we do not get exactness for all
polynomials of degree $\leq 4$. We summarize

\[
\text{\textbf{CR4}: }L(f)=\frac{1}{3}f(1/2,1/2)+\frac{1}{6}%
(f(1/2,0)+f(0,1/2)+f(1/2,1)+f(1,1/2))\text{ } 
\]
is exact for all polynomials of degree $\leq 3$. Note that $L(f)$ is also
exact for $x^{3}y$ and $xy^{3}$.

\begin{remark}
It is not possible to choose $d$ so that $\lambda =0$, since $6d^{2}-6d+2$
has no real roots.
\end{remark}

\begin{remark}
For the $n$ cube in general, $n\geq 3$, the natural extension would be to
use the center of mass of the faces of $C_{n}$ to generate $T(f)$. We leave
the details to the reader.
\end{remark}

\begin{remark}
Since $C_{n}=[0,1]\times \cdots \times [0,1]$, one could also use Simpson's
Rule as a \textbf{product rule}. However, this does \textbf{not} yield any
of the rules obtained in this section.
\end{remark}

\subsection{Connection with Interpolation}

Let

\[
p(x,y)=A_{1}x^{2}+A_{2}y^{2}+A_{3}x+A_{4}y+A_{5} 
\]
Given $0\leq a\leq 1,0\leq b\leq 1,0\leq c\leq 1,0\leq d\leq 1$, and $f(x,y)$%
, let

$P_{1}=(a,0),P_{2}=(0,b),P_{3}=(c,1),P_{4}=(1,d).$ We wish to choose the $%
A_{j}$ so that 
\[
p(P_{k})=f(P_{k}),\;k=1,...,5 
\]
where $P_{5}=(1/2,1/2)$. The corresponding coefficient matrix is $B=\left( 
\begin{array}{lllll}
a^{2} & 0 & a & 0 & 1 \\ 
0 & b^{2} & 0 & b & 1 \\ 
c^{2} & 1 & c & 1 & 1 \\ 
1 & d^{2} & 1 & d & 1 \\ 
\frac{1}{4} & \frac{1}{4} & \frac{1}{2} & \frac{1}{2} & 1
\end{array}
\right) $. If $a=1/2$, $b=1/2$, $c=1/2$, $d=1/2$, then $\det (B)=\allowbreak 
\frac{1}{16}\neq 0$. In this case there is a unique interpolant, and $%
B=\left( \allowbreak 
\begin{array}{ccccc}
\frac{1}{4} & 0 & \frac{1}{2} & 0 & 1 \\ 
0 & \frac{1}{4} & 0 & \frac{1}{2} & 1 \\ 
\frac{1}{4} & 1 & \frac{1}{2} & 1 & 1 \\ 
1 & \frac{1}{4} & 1 & \frac{1}{2} & 1 \\ 
\frac{1}{4} & \frac{1}{4} & \frac{1}{2} & \frac{1}{2} & 1
\end{array}
\right) $. Then $\left( 
\begin{array}{l}
A_{1} \\ 
A_{2} \\ 
A_{3} \\ 
A_{4} \\ 
A_{5}
\end{array}
\right) =B^{-1}\left( 
\begin{array}{l}
b_{1} \\ 
b_{2} \\ 
b_{3} \\ 
b_{4} \\ 
b_{5}
\end{array}
\right) =\allowbreak \left( 
\begin{array}{c}
2b_{2}+2b_{4}-4b_{5} \\ 
2b_{1}+2b_{3}-4b_{5} \\ 
-3b_{2}-b_{4}+4b_{5} \\ 
-3b_{1}-b_{3}+4b_{5} \\ 
b_{1}+b_{2}-b_{5}
\end{array}
\right) $, where $%
b_{1}=f(a,0),b_{2}=f(0,b),b_{3}=f(c,1),b_{4}=f(1,d),b_{5}=f(1/2,1/2)$. $%
\int_{0}^{1}\int_{0}^{1}(A_{1}x^{2}+A_{2}y^{2}+A_{3}x+A_{4}y+A_{5})dxdy=%
\frac{1}{3}(A_{1}+A_{2})+\frac{1}{2}(A_{3}+A_{4})+A_{5}=$

$\frac{1}{3}(2b_{2}+2b_{4}-4b_{5}+2b_{1}+2b_{3}-4b_{5})+\frac{1}{2}%
(-3b_{2}-b_{4}+4b_{5}-3b_{1}-b_{3}+4b_{5})+b_{1}+b_{2}-b_{5}=\allowbreak $

$\frac{1}{6}b_{2}+\frac{1}{6}b_{4}+\frac{1}{3}b_{5}+\frac{1}{6}b_{1}+\frac{1%
}{6}b_{3}=\frac{1}{3}f(1/2,1/2)+\frac{1}{6}%
(f(1/2,0)+f(0,1/2)+f(1/2,1)+f(1,1/2))$

$=L(f)$ above. Thus $L(f)=\frac{1}{3}f(1/2,1/2)+\frac{1}{6}%
(f(1/2,0)+f(0,1/2)+f(1/2,1)+f(1,1/2))$ does arise as the integral of a
unique interpolant of the form $A_{1}x^{2}+A_{2}y^{2}+A_{3}x+A_{4}y+A_{5}$.

If $a=1$, $b=0$, $c=0$, $d=1$, however, then $\det (B)=\allowbreak 0.$ Thus,
in this case, the interpolant of the form $%
A_{1}x^{2}+A_{2}y^{2}+A_{3}x+A_{4}y+A_{5}$ is not unique.

It is also natural to ask whether the cubature rule

$T(f)=\frac{1}{4}(f(1/2,0)+f(0,1/2)+f(1/2,1)+f(1,1/2))$ equals the integral
of a unique second degree interpolant to $f$ of the form 
\[
p(x,y)=A_{1}xy+A_{2}x+A_{3}y+A_{4} 
\]
We wish to choose the $A_{j}$ so that 
\[
p(P_{k})=f(P_{k}),\;k=1,...,4 
\]
The corresponding coefficient matrix is $A=\left( 
\begin{array}{llll}
0 & a & 0 & 1 \\ 
0 & 0 & b & 1 \\ 
c & 1 & c & 1 \\ 
d & d & 1 & 1
\end{array}
\right) $, and

$\det (A)=\allowbreak cab-ac-cdb-dab+dac+db$. Note that if $a=b=c=d=\frac{1}{%
2}$, then

$\det (A)=0$. Also, if $a=1$, $b=0$, $c=0$, $d=1$(gives the vertices), then
again $\det (A)=\allowbreak 0$. In either case, there is not a unique
interpolant, $p(x,y),$ in general, and thus $T(f)$ does not arise as the
integral of a unique interpolant of the form $A_{1}xy+A_{2}x+A_{3}y+A_{4}$.

\section{Polygons in the Plane}

We have already discussed cubature formulas over triangles as a special case
of the $n$ simplex, and with points on the boundary other than the vertices.
We also discussed cubature formulas over the unit square, as a special case
of the $n$ cube, and with points on the boundary other than the vertices. If 
$n>4$, and $T(f)$ is generated using the vertices, then the weighted
combination $\lambda M(f)+(1-\lambda )T(f)$ is a poor CR. We examine the
special case $n=6$.

\subsection{6 sided Regular Polygons}

Let $D$ be the regular hexagon with vertices

$(1+\sqrt{3},0),(1,1),(-1,1),(-1-\sqrt{3},0),(-1,-1),(1,-1).$ Then the
Center of Mass of

$D=(0,0)$. Using

$\int_{D}f(x,y)dA=\int_{-1}^{1}\int_{-1}^{1}f(x,y)dydx+\int_{1}^{1+\sqrt{3}%
}\int_{(x-1-\sqrt{3})/\sqrt{3}}^{-(x-1-\sqrt{3})/\sqrt{3}}f(x,y)dydx+%
\int_{-1-\sqrt{3}}^{-1}\int_{-(x+1+\sqrt{3})/\sqrt{3}}^{(x+1+\sqrt{3})/\sqrt{%
3}}f(x,y)dydx$

we list some useful integrals:

Area$(D)=\int_{D}dA=\allowbreak 4+2\sqrt{3}$

$\int_{D}xdA=\allowbreak
\int_{D}ydA=\int_{D}x^{3}dA=\int_{D}y^{3}dA=\int_{D}xydA=\int_{D}x^{2}ydA=%
\int_{D}xy^{2}dA=0$

$\int_{D}x^{2}dA=\frac{16}{3}+3\sqrt{3},\int_{D}y^{2}dA=\allowbreak \frac{4}{%
3}+\frac{1}{3}\sqrt{3}$. Let

$M(f)=(\allowbreak +2\sqrt{3})f(0,0),$

$T(f)=(4+2\sqrt{3})\frac{1}{6}(f(1+\sqrt{3},0)+f(1,1)+f(-1,1)+f(-1-\sqrt{3}%
,0)+f(-1,-1)+f(1,-1))$

$L_{\lambda }(f)=$ $\lambda M(f)+(1-\lambda )T(f)$. It follows easily that $%
L_{\lambda }$ is exact for $1,x,y$, for any $\lambda $.

$L_{\lambda }(x^{2})=\left( 4+2\sqrt{3}\right) \left( \frac{1}{6}-\frac{1}{6}%
\lambda \right) \left( \left( 1+\sqrt{3}\right) ^{2}+4+\left( -1-\sqrt{3}%
\right) ^{2}\right) $, and

$L_{\lambda }(y^{2})=4\left( 4+2\sqrt{3}\right) \left( \frac{1}{6}-\frac{1}{6%
}\lambda \right) $. To make $L_{\lambda }$ exact for $y^{2}$, say, we need

$\allowbreak 4\left( 4+2\sqrt{3}\right) \left( \frac{1}{6}-\frac{1}{6}%
\lambda \right) =\allowbreak \frac{4}{3}+\frac{1}{3}\sqrt{3}\Rightarrow
\lambda =-\frac{\frac{4}{3}+\sqrt{3}}{-\frac{8}{3}-\frac{4}{3}\sqrt{3}}%
=\allowbreak \frac{1}{4}\frac{4+3\sqrt{3}}{2+\sqrt{3}}$. But then

$\allowbreak L_{\lambda }(x^{2})=\left( 4+2\sqrt{3}\right) \left( \frac{1}{6}%
-\frac{1}{6}\lambda \right) \left( \left( 1+\sqrt{3}\right) ^{2}+4+\left( -1-%
\sqrt{3}\right) ^{2}\right) =\allowbreak $

$\frac{7}{3}\sqrt{3}+5\neq I(x^{2})=\frac{16}{3}+3\sqrt{3}$. Hence $%
L_{\lambda }$ is only exact for the class of linear polynomials.

\section{Irregular Quadrilaterals}

Let $D\subset R^{2}$ be the trapezoid with vertices $%
\{(0,0),(1,0),(0,1),(1,2)\}$(call them $P_{j}$).

$\int_{D}dA=\allowbreak \frac{3}{2}$, $\int_{D}xdA=\allowbreak \frac{5}{6}$, 
$\int_{D}ydA=\allowbreak \frac{7}{6}\Rightarrow $ Center of mass of $D$ is $%
(5/9,7/9).$

Define the following linear functionals:

\[
M(f)=\text{Area}(D)\text{ }f(5/9,7/9)=\frac{3}{2}f(5/9,7/9) 
\]
\[
T(f)=\text{Area}(D)(\frac{1}{4}\sum_{j=1}^{4}f(P_{j}))=\frac{3}{8}%
(f(0,0)+f(1,0)+f(0,1)+f(1,2) 
\]
For fixed $\lambda ,\;0\leq \lambda \leq 1$, 
\[
L_{\lambda }=\lambda M(f)+(1-\lambda )T(f) 
\]
To make $L_{\lambda }$ exact for $x$, we have $L_{\lambda }(x)=\allowbreak 
\frac{1}{12}\lambda +\frac{3}{4}=\frac{5}{6}\Rightarrow \lambda =1$. Then $%
L_{\lambda }(y)=\allowbreak \frac{1}{24}\lambda +\frac{9}{8}=\frac{7}{6}$.
However, $L_{\lambda }(xy)=\frac{35}{54}\allowbreak \neq $ $%
\int_{0}^{1}\int_{0}^{x+1}xydydx=\frac{17}{24}.$ So using the \underline{%
vertices} of $\partial (D)$ for $T(f)$ only gives exactness for \underline{%
linear} functions in general.

\begin{remark}
It is interesting to see if there are \textbf{any} weights $%
w_{1},w_{2},w_{3},w_{4},w_{5}$ such that the CR $%
L(f)=w_{1}f(5/9,7/9)+w_{2}f(0,0)+w_{3}f(1,0)+w_{4}f(0,1)+w_{5}f(1,2)$ is
exact for all polynomials of degree $\leq 2$. Above we considered the
special case where the weights would have the form $w_{1}=\frac{3}{2}$ $%
\lambda $, $w_{j}=\frac{3}{8}(1-\lambda )$, $j=2,3,4,5$. Here is the
resulting linear system: $\allowbreak w_{1}+w_{2}+w_{3}+w_{4}+w_{5}=\frac{3}{%
2}$, $\frac{5}{9}w_{1}+w_{3}+w_{5}=\frac{5}{6}$, $\frac{7}{9}%
w_{1}+w_{4}+2w_{5}=\frac{7}{6}$, $\allowbreak \frac{25}{81}w_{1}+w_{3}+w_{5}=%
\frac{7}{12}$,

$\allowbreak \frac{49}{81}w_{1}+w_{4}+4w_{5}=\frac{5}{4}$, $\allowbreak 
\frac{35}{81}w_{1}+2w_{5}=\frac{17}{24}$

The augmented matrix is $A=\left( 
\begin{array}{llllll}
1 & 1 & 1 & 1 & 1 & \frac{3}{2} \\ 
\frac{5}{9} & 0 & 1 & 0 & 1 & \frac{5}{6} \\ 
\frac{7}{9} & 0 & 0 & 1 & 2 & \frac{7}{6} \\ 
\frac{25}{81} & 0 & 1 & 0 & 1 & \frac{7}{12} \\ 
\frac{49}{81} & 0 & 0 & 1 & 4 & \frac{5}{4} \\ 
\frac{35}{81} & 0 & 0 & 0 & 2 & \frac{17}{24}
\end{array}
\right) $, and the RREF of $A$ is

$\left( 
\begin{array}{cccccc}
1 & 0 & 0 & 0 & 0 & 0 \\ 
0 & 1 & 0 & 0 & 0 & 0 \\ 
0 & 0 & 1 & 0 & 0 & 0 \\ 
0 & 0 & 0 & 1 & 0 & 0 \\ 
0 & 0 & 0 & 0 & 1 & 0 \\ 
0 & 0 & 0 & 0 & 0 & 1
\end{array}
\right) $. Hence the system has no solution, and thus there are \textbf{no}
weights $w_{1},w_{2},w_{3},w_{4},w_{5}$ such that the CR $%
L(f)=w_{1}f(5/9,7/9)+w_{2}f(0,0)+w_{3}f(1,0)+w_{4}f(0,1)+w_{5}f(1,2)$ is
exact for all

polynomials of degree $\leq 2$. One can, however, choose weights so that $%
L(f)$ is exact for all polynomials of degree $\leq 2$, except for $xy$, say.
The corresponding augmented matrix is

$\left( 
\begin{array}{llllll}
1 & 1 & 1 & 1 & 1 & \frac{3}{2} \\ 
\frac{5}{9} & 0 & 1 & 0 & 1 & \frac{5}{6} \\ 
\frac{7}{9} & 0 & 0 & 1 & 2 & \frac{7}{6} \\ 
\frac{25}{81} & 0 & 1 & 0 & 1 & \frac{7}{12} \\ 
\frac{49}{81} & 0 & 0 & 1 & 4 & \frac{5}{4}
\end{array}
\right) $, which has row echelon form: $\left( 
\begin{array}{cccccc}
1 & 0 & 0 & 0 & 0 & \frac{81}{80} \\ 
0 & 1 & 0 & 0 & 0 & \frac{23}{240} \\ 
0 & 0 & 1 & 0 & 0 & \frac{17}{120} \\ 
0 & 0 & 0 & 1 & 0 & \frac{29}{240} \\ 
0 & 0 & 0 & 0 & 1 & \frac{31}{240}
\end{array}
\right) $
\end{remark}

\subsubsection{Other points on Boundary of Trapezoid}

We will try using other points on $\partial (D)$ to generate $T(f)-$say $%
(a,0),(0,b),(1,c),(d,d+1),$ with $0\leq a\leq 1,0\leq b\leq 1,0\leq c\leq
2,0\leq d\leq 1$. Then 
\[
L_{\lambda }(f)=\lambda \frac{3}{2}f(5/9,7/9)+(1-\lambda )\frac{3}{8}%
(f(a,0)+f(1,c)+f(0,b)+f(d,d+1) 
\]

\begin{itemize}
\item  $L_{\lambda }(x)=\allowbreak \frac{5}{6}\lambda +\frac{3}{8}\left(
1-\lambda \right) \left( a+1+d\right) $

\item  $L_{\lambda }(y)=\allowbreak \frac{7}{6}\lambda +\frac{3}{8}\left(
1-\lambda \right) \left( c+b+d+1\right) $

\item  $L_{\lambda }(xy)=\allowbreak \frac{35}{54}\lambda +\frac{3}{8}\left(
1-\lambda \right) \left( c+d\left( d+1\right) \right) $

\item  $L_{\lambda }(x^{2})=\allowbreak \frac{25}{54}\lambda +\frac{3}{8}%
\left( 1-\lambda \right) \left( a^{2}+1+d^{2}\right) $

\item  $L_{\lambda }(y^{2})=\allowbreak \frac{49}{54}\lambda +\frac{3}{8}%
\left( 1-\lambda \right) \left( c^{2}+b^{2}+\left( d+1\right) ^{2}\right) $
\end{itemize}

We set each of the above expressions in $\lambda ,a,b,c,d$ equal to the
corresponding integral of $f$ over $D$ to yield the system:

$\frac{5}{6}\lambda +\frac{3}{8}\left( 1-\lambda \right) \left( a+1+d\right)
=\frac{5}{6},\frac{7}{6}\lambda +\frac{3}{8}\left( 1-\lambda \right) \left(
c+b+d+1\right) =\frac{7}{6}$

$\frac{35}{54}\lambda +\frac{3}{8}\left( 1-\lambda \right) \left( c+d\left(
d+1\right) \right) =\frac{17}{24},\frac{25}{54}\lambda +\frac{3}{8}\left(
1-\lambda \right) \left( a^{2}+1+d^{2}\right) =\frac{7}{12}$

$\frac{49}{54}\lambda +\frac{3}{8}\left( 1-\lambda \right) \left(
c^{2}+b^{2}+\left( d+1\right) ^{2}\right) =\frac{5}{4}.$

We used Maple to find the following Grobner basis:

\[
\{392\lambda -163,9a+9d-11,81b-99d+20,180d-191+81c,-22671d+6583+18549d^{2}\} 
\]
The last equation has solutions $d=\frac{11}{18}\pm \frac{1}{458}\sqrt{3893}%
\approx .\,74734,\allowbreak .\,47488$

We shall use $d=\frac{11}{18}+\frac{1}{458}\sqrt{3893}$. Solving the other
equations yields

$180d-191+81c=0\Rightarrow c=1-\frac{10}{2061}\sqrt{3893}\approx \allowbreak
.\,69726$

$81b-99d+20=0\Rightarrow b=\frac{1}{2}+\frac{11}{4122}\sqrt{3893}\approx
\allowbreak .\,6665$

$9a+9d-11=0\Rightarrow a=\frac{11}{18}-\frac{1}{458}\sqrt{3893}\approx
\allowbreak .\,47488$

$392\lambda -163=0\Rightarrow \lambda =\frac{163}{392}\approx \allowbreak
.\,41582$. We summarize

\textbf{CR5}: Let $D\subset R^{2}$ be the trapezoid with vertices $%
\{(0,0),(1,0),(0,1),(1,2)\}$.

Let $a=\frac{11}{18}-\frac{1}{458}\sqrt{3893}\approx \allowbreak .\,47488$, $%
b=\frac{1}{2}+\frac{11}{4122}\sqrt{3893}\approx \allowbreak .\,6665$,

$c=1-\frac{10}{2061}\sqrt{3893}\approx \allowbreak .\,69726$, $d=\frac{11}{18%
}+\frac{1}{458}\sqrt{3893}\approx \allowbreak .\,74734$, and $\lambda =\frac{%
163}{392}\approx \allowbreak .\,41582$. Then $L(f)=\lambda \frac{3}{2}%
f(5/9,7/9)+(1-\lambda )\frac{3}{8}(f(a,0)+f(1,c)+f(0,b)+f(d,d+1))$ is exact
for all polys. of degree $\leq 2$.

Note that $L(x^{3})=\frac{336001}{762048}\approx .\,44092$, while $%
I(x^{3})=\allowbreak \frac{9}{20}=\allowbreak .\,45$. So $L$ is not exact
for cubics.

\section{Unit Disc}

Let $D$ be the unit disc. The $n$ roots of unity $z_{k}=e^{2\pi
ik/n},k=1,...,n,$ can be used to generate $T(f)=\frac{\pi }{n}%
\sum_{k=1}^{n}f(\cos ($ $2\pi k/n,\sin (2\pi k/n))$. Since the center of
mass of $D=(0,0)$ and area$(D)=\pi $, the analogy of our formulas for the
triangle and the square is $M(f)=\pi f(0,0)$, $L_{\lambda }(f)=$ $\lambda
\pi f(0,0)+(1-\lambda )\frac{\pi }{n}\sum_{k=1}^{n}f(\cos ($ $2\pi k/n,\sin
(2\pi k/n)).$ It is not hard to show, however, that the only good choice is $%
n=4$. This gives the formulas 
\[
L_{\lambda }(f)=\lambda \pi f(0,0)+(1-\lambda )\frac{\pi }{4}%
(f(1,0)+f(0,1)+f(-1,0)+f(0,-1)) 
\]
We list the following integrals without proof:

\begin{itemize}
\item  If $m$ and $n$ are even whole numbers, then 
\[
\int \int_{D}x^{m}y^{n}dA=\frac{\pi }{m+n+2}\left( \frac{(n-1)!(m-1)!}{%
2^{m+n-3}(\frac{n}{2}-1)!(\frac{m}{2}-1)!(\frac{m+n}{2})!}\right) 
\]

\item  If $m$ is an even whole number and $n=0$, then
\end{itemize}

\[
\int \int_{D}x^{m}y^{n}dA=\frac{\pi }{m+2}\left( \frac{(m-1)!}{2^{m-2}(\frac{%
m}{2}-1)!(\frac{m}{2})!}\right) 
\]

\begin{itemize}
\item  If $n$ is an even whole number and $m=0$, then
\end{itemize}

\[
\int \int_{D}x^{m}y^{n}dA=\frac{\pi }{n+2}\left( \frac{(n-1)!}{2^{n-2}(\frac{%
n}{2}-1)!(\frac{n}{2})!}\right) 
\]

\begin{itemize}
\item  If $m$ and/or $n$ is an odd whole number, then
\end{itemize}

\[
\int \int_{D}x^{m}y^{n}dA=0 
\]

$\allowbreak $It is then easy to prove that the choice $\lambda =\frac{1}{2}$
gives

\[
\text{\textbf{CR6:}}\mathbf{\ }L(f)=\frac{\pi }{2}f(0,0)+\frac{\pi }{8}%
(f(1,0)+f(0,1)+f(-1,0)+f(0,-1))\text{ } 
\]
is exact for all polynomials of degree $\leq 3.$

\section{Summary of Cubature Rules}

These are the CRs we derived as a generalization of Simpson's Rule for
functions of one variable. If $D_{n}$ is a polygonal region in $R^{n}$, let $%
P_{0},...,P_{n}$ denote the $n+1$ vertices of $D_{n}$ and $P_{n+1}=$ Center
of Mass of $D_{n}$. The generalization is based on the weighted combination $%
L_{\lambda }=\lambda M(f)+(1-\lambda )T(f)$, where $M(f)=$ Vol$(D_{n})$ $%
f(P_{n+1}),$ $T(f)=Vol(D_{n})(\frac{1}{m+1}\sum_{j=0}^{m}f(P_{j}))$.

Unless noted otherwise, all of the rules have the following property:

All of the weights are positive, all of the knots lie inside the region, and
all but one of the knots lies on the boundary of the region. It is desirable
to have as many knots as possible on $\partial (D_{n})$ if one subdivides
the region and compounds the CR.

\subsection{$n$ Simplex $T_{n}$}

\begin{itemize}
\item  Using the vertices, $L(f)=\frac{n+1}{(n+2)n!}f(P_{n+1})+\frac{1}{%
(n+2)!}\sum_{j=0}^{n}f(P_{j})$ is exact for all polynomials of degree $\leq
2 $.

\item  Using points other than the vertices: Letting $\{Q_{k}\}$ denote the
center of mass of the faces of $T_{n},L(f)=-\left( n-2\right) \frac{n+1}{n+2}%
\frac{1}{n!}f(1/(n+1),...,1/(n+1))+\frac{n^{2}}{(n+2)!}%
\sum_{k=1}^{n+1}f(Q_{k})$ is exact for all polynomials of degree $\leq 2$.
All but one weight is positive if $n>2$.
\end{itemize}

\subsection{Unit $n$ Cube $C_{n}$}

\begin{itemize}
\item  General $n$: Let $P_{0},...,P_{m-1}$ denote the $m$ vertices of $%
C_{n},m=2^{n}$. Then

$L(f)=\frac{2}{3}f(1/2,...,1/2)+\frac{1}{3}\frac{1}{2^{n}}%
\sum_{j=0}^{2^{n}-1}f(P_{j})$ is exact for all polynomials of degree $\leq 3$%
.

\item  $n=2:$ $L(f)=\frac{1}{3}f(1/2,1/2)+\frac{1}{6}%
(f(1/2,0)+f(0,1/2)+f(1/2,1)+f(1,1/2))$ is exact for all polynomials of
degree $\leq 3$.
\end{itemize}

\subsection{Unit Disc $D$}

$L(f)=\frac{\pi }{2}f(0,0)+\frac{\pi }{8}(f(1,0)+f(0,1)+f(-1,0)+f(0,-1))$ is
exact for all polynomials of degree $\leq 3$.

\subsection{A Trapezoid in the Plane}

Let $D\subset R^{2}$ be the trapezoid with vertices $%
\{(0,0),(1,0),(0,1),(1,2)\}$.

Let $a=\frac{11}{18}-\frac{1}{458}\sqrt{3893}\approx \allowbreak .\,47488$, $%
b=\frac{1}{2}+\frac{11}{4122}\sqrt{3893}\approx \allowbreak .\,6665$,

$c=1-\frac{10}{2061}\sqrt{3893}\approx \allowbreak .\,69726$, $d=\frac{11}{18%
}+\frac{1}{458}\sqrt{3893}\approx \allowbreak .\,74734$, and $\lambda =\frac{%
163}{392}\approx \allowbreak .\,41582$. Then $L(f)=\lambda \frac{3}{2}%
f(5/9,7/9)+(1-\lambda )\frac{3}{8}(f(a,0)+f(1,c)+f(0,b)+f(d,d+1))$ is exact
for all polys. of degree $\leq 2$.

\end{document}